\documentclass{article}
\usepackage[utf8]{inputenc}
\usepackage{amsfonts}
\usepackage{epigraph}

\usepackage{amsfonts}
\usepackage{amssymb}
\usepackage{amsmath}
\usepackage{amssymb}
\usepackage{float}
\usepackage{amsthm}
\usepackage[pdftex]{graphicx}

\usepackage{fullpage}
\usepackage{wrapfig}

\newtheorem{proposition}{Proposition}

\author{D. Cherkashin$^{a,c}$, F. Petrov$^{a,b}$, V. Sokolov$^b$}

\title{On a question of Sidorenko}

\begin{document}
\maketitle
{\let\thefootnote\relax\footnote{
a) St. Petersburg State University;
b) St. Petersburg Department of
V.~A.~Steklov Institute of Mathematics of
the Russian Academy of Sciences;
c) St. Petersburg Branch of Higher School of Economics.
E-mails: matelk@gmail.com,  f.v.petrov@spbu.ru,visoksok@gmail.com}}

In a very recent paper~\cite{sidorenko2018tur} Sidorenko stated the following problem:
\begin{quote}
Let $G_k$ be a graph whose vertices are functions $f : \mathbb{Z}_k \to \mathbb{Z}_k$. A pair of
vertices $\{f, g\}$ forms an edge in $G_k$ if $f - g$ is a bijection. 
Lemma 2 restates the
fact that $G_k$ has no triangles when $k$ is even. For odd $k$, the problem of counting
triangles in $G_k$ has been solved asymptotically in~\cite{eberhard2015additive}. Let $p(k)$ be the smallest
prime factor of $k$. The $p(k)$ functions $f_0, f_1, \dots, f_{p(k)-1}$, where $f_i(j) := i \cdot j \mbox{ mod } k$, form a complete subgraph in $G_k$. It is very tempting to conjecture that $p(k)$ is indeed the size of the largest clique in $G_k$. We know that this is
true for even $k$ and for prime $k$. Computer search confirms that this is also true for $k = 9$.
\end{quote}
It turns out that there is a counterexample for $k = 15$.
\[
\begin{array}{ccccccccccccccc}
0 & 0 & 0 & 0 & 0 & 0 & 0 & 0 & 0 & 0 & 0 & 0 & 0 & 0 & 0 \\
0 & 1 & 2 & 3 & 4 & 5 & 6 & 7 & 8 & 9 & 10 & 11 & 12 & 13 & 14 \\
0 & 9 & 3 & 2 & 13 & 11 & 10 & 12 & 4 & 6 & 8 & 14 & 7 & 5 & 1 \\
0 & 12 & 4 & 11 & 10 & 9 & 5 & 2 & 6 & 14 & 7 & 3 & 13 & 1 & 8
\end{array}
\]

It is found by computer search,
we do not see any specific structure in it.

More intensive computer 
search allowed to get the examples of four functions
for $k=21$ and $k=27$. 
Two of four functions are again $f_0(x)=0$
and $f_1(x)=x$,
$0\leqslant 
x\leqslant k-1$,
(that may
be always assumed a priori),  and two other are, starting from the
value at 0:
$$
\begin{array}{ccccccccccccccccccccc}
    13& 11& 14& 0& 2& 1& 5& 7& 3& 10& 15& 17& 16& 20& 4& 18& 9& 19& 12& 6& 8\\ 
14& 5& 4& 13& 9& 18& 
2& 15& 6& 10& 17& 1& 11& 19& 
8& 3& 7& 12& 0& 16& 20
\end{array}
$$
for $k=21$ and 
$$
\begin{array}{ccccccccccccccccccccccccccc}
  12& 17& 11& 20& 5& 19& 1& 9& 0& 13& 15& 18& 6& 10& 22& 3& 2& 8& 14& 25& 4& 24& 21& 16& 7& 23& 26\\
4& 6& 5& 15& 19& 18& 3& 13& 24& 16& 20& 1& 7& 0& 8& 11& 9& 17& 26& 21& 2& 12& 14& 22& 25& 23& 10
\end{array}
$$
for $k=27$.

Denote by $\omega (k)$ the size of the largest clique in $G_k$. We have the following general
\begin{proposition}
\[
\omega(nm) \geq \min (\omega(n), \omega (m)).
\]
\end{proposition}
\begin{proof}
Consider arbitrary corresponding functions $f_1, \dots, f_{\omega(n)}$, $g_1, \dots g_{\omega(m)}$.
For $t$ from 1 to $\min (\omega(n), \omega (m))$ put
\[
q_t (im+j) := f_t(i) m + g_t(j).
\]
Suppose that $q_{t_1} - q_{t_2}$ is not a bijection. Then for some different $a$, $b$ we have 
\[
q_{t_1}(a)  - q_{t_2} (a) = q_{t_1}(b) - q_{t_2}(b) \mbox{ mod } nm.
\]
Let $a = im + j$, $b = rm + s$. Then
\begin{equation}
    f_{t_1}(i) m+ g_{t_1}(j) - f_{t_2}(i) m - g_{t_2}(j) = f_{t_1}(r) m+ g_{t_1}(s) - f_{t_2}(r) m - g_{t_2}(s) \mbox{ mod } nm.
    \label{1}
\end{equation}
Modulo $m$ we have
\[
g_{t_1}(j) - g_{t_2}(j) = g_{t_1}(s) - g_{t_2}(s) \mbox{ mod } m,
\]
so $j=s$, because $g_{t_1}$ and $g_{t_2}$ were connected in $G_m$.
Dividing~\eqref{1} by $m$ we get
\[
f_{t_1}(i) - f_{t_2}(i) = f_{t_1}(r) - f_{t_2}(r) \mbox{ mod } n.
\]
Hence $i = r$, and $a = b$, contradiction. So $\{q_t\}$ is a clique of size $\min (\omega(n), \omega (m))$.
\end{proof}

In particular, this gives a lower estimate 
$\omega(15^{k_1}21^{k_2}27^{k_3} n)\geqslant 4$
for any $n$ coprime to 
6 and integers $k_1,k_2,k_3\geqslant 0$.
We expect that $\omega(3n)\geqslant 4$ for any
odd $n>3$.

\bibliographystyle{plain}
\bibliography{main}

\begin{thebibliography}{1}

\bibitem{eberhard2015additive}
Sean Eberhard, Freddie Manners, and Rudi Mrazovi{\'c}.
\newblock Additive triples of bijections, or the toroidal semiqueens problem.
\newblock {\em arXiv preprint arXiv:1510.05987}, 2015.

\bibitem{sidorenko2018tur}
Alexander Sidorenko.
\newblock On {T}ur{\'a}n problems for {C}artesian products of graphs.
\newblock {\em arXiv preprint arXiv:1812.01581}, 2018.

\end{thebibliography}

\end{document}